\documentclass[reqno]{amsart}
\usepackage{amscd,amsfonts,amssymb}
\usepackage{graphicx}
\usepackage{color}

\numberwithin{equation}{section}
\newtheorem{Theorem}{Theorem}[section]
\newtheorem{Lemma}{Lemma}[section]

\theoremstyle{definition}

\theoremstyle{remark}
\newtheorem{Remark}{Remark}[section]

\renewcommand{\r}{\rho}

\renewcommand{\i}{\varepsilon}

\renewcommand{\u}{{\bf u}}
\renewcommand{\H}{{\bf H}}
\newcommand{\R}{{\mathbb R}}
\newcommand{\Dv}{{\rm div}}

\newcommand{\na}{\nabla}
\newcommand{\dl}{\delta}

\def\f{\frac}
\renewcommand{\O}{\Omega}
\def\ov{\overline}

\def\hf1{^\f{1}{1-\xi^2}}

\def\be{\begin{equation}}
\def\en{\end{equation}}
\def\bs{\begin{split}}
\def\es{\end{split}}

\author{Xianpeng Hu \and Dehua Wang}
\address{Department of Mathematics, University of Pittsburgh,
                           Pittsburgh, PA 15260, USA.}
\email{xih15@pitt.edu}
\address{Department of Mathematics, University of Pittsburgh,
                           Pittsburgh, PA 15260, USA.}
\email{dwang@math.pitt.edu}

\title[Low Mach Number Limit]
{Low Mach Number Limit of Viscous Compressible Magnetohydrodynamic
Flows} 
\date{\today}

\begin{document}

\begin{abstract}
The relationship between the  compressible
magnetohydrodynamic flows with low Mach number and the incompressible magnetohydrodynamic flows is investigated.
More precisely,  the convergence of weak solutions of the
compressible isentropic viscous magnetohydrodynamic equations to the weak
solutions of the incompressible viscous magnetohydrodynamic equations is proved as the
density becomes constant and the Mach number goes to zero, that is, the corresponding incompressible limits are justified when the  spatial domain is a periodic domain, the whole space, or a bounded domain.
\end{abstract}

\maketitle

\section{Introduction}

Studies on magnetohydrodynamic flows always involve a choice at the
onset to describe the system entirely in the context of either
incompressible magnetohydrodynamics (MHD), or compressible MHD. For
example,  theoretic studies on turbulence have a particular leaning
toward  the incompressible model. This preference has largely
been based on the benefits and advantages of the similarity of
incompressible MHD to its hydrodynamic counterparts, and the
practical consideration of limited computational resources.
However, when the density of a flow is no longer invariant, the
flow become much more complicated not only from the physical
viewpoint, but also from the mathematical consideration, see
\cite{Ca, hw1, hw2, KL, LL} and references therein. Thus, it is a
natural problem to consider the relation between the
incompressible MHD and the compressible MHD. The equations of the
isentropic compressible viscous magnetohydrodynamic flows in $N$
spatial dimensions have the following form (\cite{Ca, KL, LL}):
\begin{equation}\label{1}
\begin{cases}
{\widetilde{\r}}_t+\Dv(\widetilde{\r}\,\widetilde{\u})=0,\\
(\widetilde{\r}\,\widetilde{\u})_t
   +\Dv\left(\widetilde{\r}\,\widetilde{\u}\otimes\widetilde{\u}\right)
   +\nabla \widetilde{p}(\widetilde{\r}) =(\na \times \widetilde{\H})\times
\widetilde{\H}+\widetilde{\mu} \Delta
\widetilde{\u}+\widetilde{\lambda}\nabla\Dv\widetilde{\u},\\
\widetilde{\H}_t-\nabla\times(\widetilde{\u}\times\widetilde{\H})
   =-\nabla\times(\widetilde{\nu}\,\nabla\times\widetilde{\H}),\quad \Dv\widetilde{\H}=0,
\end{cases}
\end{equation}
where $\widetilde{\mu}>0$ is the  shear viscosity,  $\widetilde{\lambda}$ is the
bulk viscosity  satisfying  $2\widetilde{\mu}+N\widetilde{\lambda}>0$ , $\widetilde{\nu}>0$ is
the magnetic viscosity;  and $\widetilde{\r}$ denotes the density,
$\widetilde{\u}\in\R^N$ the velocity, $\widetilde{\H}\in\R^N$ the
magnetic field,
$\widetilde{p}(\widetilde{\r})=a\widetilde{\r}^\gamma$ the
pressure with constant $a>0$ and the adiabatic exponent
$\gamma>1$. The symbol $\otimes$ denotes the Kronecker tensor
product. The first equation in \eqref{1} is called the continuity
equation and the third equation in \eqref{1} is called the
induction equation.

From the physics point of view, the compressible flow behaves asymptotically like an
incompressible flow when the density is almost constant, and the
velocity and the magnetic field are small, in a large time
scale. More precisely, we scale $\widetilde{\r}$, $\widetilde{\u}$,
and $\widetilde{\H}$ in the following way:
\begin{equation}\label{2}
\widetilde{\r}=\r(x,\i t),\quad \widetilde{\u}=\i\u(x,\i t),\quad
\widetilde{\H}=\i\H(x,\i t),
\end{equation}
and we assume that the coefficients $\widetilde{\mu}$,
$\widetilde{\lambda}$, and $\widetilde{\nu}$ are small and scaled as:
\begin{equation}\label{3}
\widetilde\mu=\i\mu_{\i}, \quad
\widetilde\lambda=\i\lambda_{\i},\quad \widetilde\nu=\i\nu_{\i},
\end{equation}
where $\i\in(0,1)$ is a small parameter and the normalized
coefficients $\mu_{\i}$, $\lambda_{\i}$ and $\nu_{\i}$ satisfy
\begin{equation}\label{4}
\mu_{\i}\rightarrow \mu,\quad  \lambda_{\i}\rightarrow \lambda,\quad
\nu_{\i}\rightarrow \nu, \textrm{ as } \i \to 0+,
\end{equation}
with $\mu>0,\, 2\mu+N\lambda>0,\, \mathrm{and}\, \nu>0$. Such a
scaling as \eqref{3} ensures that the limit sysyem as
$\i\rightarrow 0$ is not of an Euler type. Also notice that
the parameter $\i$ in the front of the magnetic field $\H$ in
\eqref{2} can be understood as the reciprocal of Alf\'{v}en
number(\cite{MB}). Under those scalings, system \eqref{1} yields
\begin{equation}\label{c}
\begin{cases}
{\r}_t+\Dv(\r\u)=0,\\
(\r\u)_t+\Dv(\r\u\otimes\u)-\mu_{\i}\Delta\u-\lambda_{\i}\nabla\Dv\u+\f{a}{\i^2}\nabla\r^{\gamma}=
(\nabla\times\H)\times\H,\\
\H_t-\nabla\times(\u\times\H)=-\nabla\times(\nu_{\i}\nabla\times\H),\quad
\Dv\H=0.
\end{cases}
\end{equation}
The existence of global weak
solutions to \eqref{c} has been investigated in Hu-Wang \cite{hw1} (and in Hu-Wang \cite{hw2} for the non-isentropic case). From the
mathematical point of view, it is reasonable to expect that, as $\r \to 1$,
the first equation in \eqref{c} yields the limit:
$\Dv\u=0$, which is the incompressible condition of a fluid, and the
 first two terms in the second equation of \eqref{c}  become
 $$\u_t+\Dv(\u\otimes\u)=\u_t+(\u\cdot\nabla)\u.$$
  On the other hand, the incompressible MHD equations read
\begin{equation}\label{in}
\begin{cases}
\u_t+(\u\cdot\nabla)\u-\mu\Delta\u+\nabla p=(\nabla\times\H)\times\H,\\
\H_t-\nabla\times(\u\times\H)=-\nabla\times(\nu\nabla\times\H),\\
\Dv\u=0,\quad \Dv\H=0.
\end{cases}
\end{equation}
Thus, roughly speaking, it is also reasonable to expect from the
mathematical point of view that weak solutions of \eqref{c}
converge in certain suitable functional spaces to the weak
solutions of \eqref{in} as $\r$ goes to a constant such as $1$ and
$\i$ goes to 0, and the hydrostatic pressure $p$ in
\eqref{in} is the ``limit" of $(\r^\gamma-1)/\i^2$ in
\eqref{c}.
This paper is devoted to the rigorous justification of the convergence of
that incompressible limit (i.e., the low Mach number limit)
for global weak solutions of the compressible isentropic MHD equations.

In this paper, we  shall establish the incompressible limit of \eqref{c}
in three types of spatial domains:
the torus $\mathbb{T}$ (in this
case, all the functions are defined on $\R^N$ and assumed to be
periodic with period $2\pi$ for all directions, that is,
$\mathbb{T}=[0,2\pi]^N$), the whole space $\R^N$,
and  a sufficiently smooth bounded domain $\O\in\R^N$, $N=2, 3$.
The study in the bounded smooth domain
with no-slip boundary condition on the velocity is much harder
than that in other two cases, because in bounded domains, there are
extra difficulties arising from the appearance of the boundary
layers, and the subtle interactions between dissipative effects and
wave propagation near the boundary, and hence requires a different approach.
We remark that the incompressible limits for compressible
isentropic Navier-Stokes equations have been investigated in
\cite{LM} for the whole space $\R^N$ and the periodic domain using the group
method,  and in \cite{DGLM} for a bounded domain.
 These results have been extended by others, such as \cite{BDGL,DT,
DG,  MN, WJ}. We also notice that in \cite{DH}, convergence
results were proved for well-prepared data as long as the solution of
incompressible limit  is suitably smooth. For the case of non-isentropic flows,
see \cite{FN1, FN2} for some recent studies.  For other related studies on the incompressible limits of viscous and inviscid flows,
see \cite{Al, Dan, Ebin, Hoff2, KM1, KM2, Lin, MS1, MS2, Sch1, Sch2}
and the references in \cite{FN1}. Comparing with
those works on the compressible Navier-Stokes equations, we will
encounter extra difficulties in studying the compressible MHD
equations. More precisely, besides the possible oscillation of the
density, the appearance of the boundary layer and the interactions
between dissipative effects and wave propagation, the appearance of
the magnetic field and the coupling effect between the hydrodynamic
motion and the magnetic field should also been taken into
considerations with new estimates.
We will overcome all these difficulties by using the group method,
Strichartz's estimate, and the weak convergence method to
establish the convergence of weak solutions of the compressible
isentropic MHD equations \eqref{c} to weak solutions of the
incompressible MHD equations \eqref{in} as the density goes to a
constant and $\i$ goes to $0$ in the periodic case and the whole
space case. More precisely,  we will show that, for any fixed
$T>0$, in the periodic case, the incompressible part of the
velocity strongly converges to a divergence-free vector field in
$L^2([0,T], L^2(\mathbb{T}))$ while the gradient part of the
velocity converges weakly to zero; and in the whole space case,
due to Strichartz's estimate, the gradient part of the velocity
converges strongly to $0$ in $L^2([0,T],
L_{loc}^2(\mathbb{\R^N}))$, while the strong convergence of the
incompressible part of the velocity only holds in the local sense.
However, this method does not apply to the case of bounded domains
because of subtle interactions between dissipative effects and
wave propagation near the boundary. Instead, we will use the
spectral analysis of the semigroup generated by the dissipative
wave operator, together with Duhamel's principle. Finally, we
remark that the incompressible flow also can be derived from the
vanishing Debye length type limit of a compressible flows with a
Poisson damping. We refer the interested readers to \cite{mm1,mm2}.

We organize the rest of the paper as follows. In Section 2, we will
give the setting of our problem and state our main results. In
Section 3, we discuss the convergence of the incompressible limit in
the periodic case. In Section 4, we will investigate the convergence
of the incompressible limit in the whole space $\R^N$. Finally, in
Section 5, we will study the convergence of the incompressible limit
in the bounded domain.

\bigskip

\section{Main Results}

In this section, we  describe the setting of our
problem and state our main results. First, we denote by $P$ the
orthogonal projection onto incompressible vector fields, i.e.
$$v=Pv+Qv, \quad\text{with}\quad \Dv(Pv)=0, \; \textrm{curl}(Qv)=0,$$ for all $v\in L^2$. 
Indeed, in view of results in \cite{GG}, we
know that the operators $P$ and $Q$ are linear bounded operators in
$W^{s,p}$ for all $s\ge 0$ and $1<p<\infty$ in the whole space or
bounded domains with smooth boundaries. Second, let us explain the
notation of weak solutions to the incompressible MHD equations as
follows: Given the initial conditions $\u_0\in L^2$, $\H_0\in L^2$
such that $\Dv\u_0=0$ and $\Dv\H_0=0$, $(\u,\H)$ is a weak solution
of \eqref{in} satisfying
\begin{equation}\label{9}
\u|_{t=0}=\u_0, \quad \H|_{t=0}=\H_0,
\end{equation}
where
$$\u\in C([0,T]; L^2_{weak})\cap L^2([0,T]; H^1(\O)),\;
\H\in C([0,T]; L^2_{weak})\cap L^2([0,T]; H^1(\O)),$$
if for all $T<\infty$,
 $\psi\in C^\infty_0(\O)$ with $\Dv\psi=0$, and
$\varphi\in C_0^\infty([0,T))$, we have
\begin{equation*}
\begin{split}
&\psi(0)\int_\O\u_0\varphi dx
+\int_0^t \psi'(t)\int_\O \u\cdot\varphi dxdt
+\int_0^t \psi(t)\int_\O\left(\u_i\partial_i\varphi_j\u_j
-\mu \nabla\u:\nabla\varphi\right) dxdt \\
&=-\int_0^t\int_\O\psi(\nabla\times\H)\times\H\cdot\varphi dxdt,
\end{split}
\end{equation*}
and
\begin{equation*}
\begin{split}
&\psi(0)\int_\O\H_0\varphi dx
+\int_0^t \psi'(t)\int_\O\H\cdot\varphi dxdt
+\int_0^t\psi(t)\int_\O(\u\times\H)\cdot(\nabla\times\varphi) dxdt\\
&=\nu\int_0^t\psi(t)\int_\O(\nabla\times\H)\cdot(\nabla\times\varphi)dxdt.
\end{split}
\end{equation*}
For more details as to the existence and  regularity of weak
solutions to the incompressible MHD equations, we refer the readers to
\cite{DL, ST}. Now, we can state our main results case by case.

\subsection{The periodic case}

Let us begin with the periodic case. We consider a sequence of
global weak solutions $(\r_{\i}, \u_{\i}, \H_{\i})$ of the compressible
MHD equations \eqref{c} in $\mathbb{T}$ and  assume that
$$\r_{\i}\in L^{\infty}([0,T]; L^\gamma(\mathbb{T})), \;
  \u_{\i}\in L^2([0,T]; H^1(\mathbb{T})),$$
$$\r_{\i}|\u_{\i}|^2 \in L^{\infty}([0,T]; L^1(\mathbb{T})), \;
\r_{\i}\u_{\i}\in C\left([0,T]; L_{weak}^{\f{2\gamma}{\gamma+1}}\right), $$
$$\H_{\i}\in L^{2}([0,T]; H^{1}(\mathbb{T}))
  \cap C([0,T]; L^2_{weak}(\mathbb{T})),$$
for all $T\in(0,\infty)$, where $C([0,T]; L^p_{weak})$ denotes the functions
which are continuous with respect to $t\in[0,T]$ with values in
$L^p$ endowed with the weak topology. We require \eqref{c} to hold
in the sense of distributions. Finally, we prescribe initial
conditions
\begin{equation}\label{6}
\r_{\i}|_{t=0}=\r^0_{\i}, \quad
\r_{\i}\u_{\i}|_{t=0}=m_{\i}^0=\r_{\i}^0 \u_{\i}^0,\quad
 \H_{\i}|_{t=0}=\H_{\i}^0,
\end{equation}
where $\r_{\i}^0\ge 0$, $\r_{\i}^0\in L^\gamma(\mathbb{T})$,
$m_{\i}^0\in L^{2\gamma/(\gamma+1)}(\mathbb{T})$, $m_{\i}^0=0$ on
$\{ \r_{\i}^0=0\}$, $\r_{\i}^0|\u_{\i}^0|^2\in L^1(\mathbb{T})$, and
$\H_{\i}^0\in L^2(\mathbb{T})$. Furthermore, we assume that
$\sqrt{\r_{\i}^0}\u_{\i}^0$ and $\H_\i^0$ converge weakly in $L^2$
to $\u_0$ and $\H_0$ respectively, and that we have
\begin{equation}\label{7}
\begin{split}
&\f{1}{2}\int_{\mathbb{T}}\left(\r_{\i}^0|\u_{\i}^0|^2+|\H_{\i}^0|^2\right)dx
+\f{a}{\i^2(\gamma-1)}\int_{\mathbb{T}}\left((\r_{\i}^0)^{\gamma}
-\gamma\r_{\i}^0(\ov{\r_{\i}})^{\gamma-1}+(\gamma-1)(\ov{\r_{\i}})^{\gamma}\right)
\le C,\\
&\ov{\r_{\i}}=(2\pi)^{-N}\int_{\mathbb{T}}\r_{\i}^0 dx \rightarrow 1, \quad\text{as }
\i\to 0,
\end{split}
\end{equation}
where and hereafter $C$ denotes a generic positive constant
independent of $\i$. Notice that \eqref{7} implies that,
roughly speaking, $\r_{\i}^0$ is of order $\ov{\r_{\i}}+O(\i)$.
We assume finally that the total energy is conserved in the sense:
\begin{equation}\label{8}
E_{\i}(t)+\int_0^t\!\!\!\!D_{\i}(s)ds\le
E^0_{\i},\quad\textrm{a.e}\quad t\in[0,T],
\end{equation}
where
$$E_{\i}=\f{1}{2}\int_\O\left(\r_{\i}|\u_{\i}|^2+|\H_{\i}|^2+\f{a}{\i^2(\gamma-1)}\r_{\i}^\gamma\right) dx,$$
$$D_{\i}=\int_\O\left(\mu_{\i}|D\u_{\i}|^2+\lambda_{\i}(\Dv\u_{\i})^2+\nu_{\i}|\nabla\times\H_{\i}|^2\right) dx,$$
and
$$E_{\i}^0=\f{1}{2}\int_\O\left(\r^0_{\i}|\u^0_{\i}|^2+|\H^0_{\i}|^2+\f{a}
{\i^2(\gamma-1)}(\r^0_{\i})^\gamma\right) dx,$$
where $\O$ is equal to $\mathbb{T}$ in the periodic case, and later is the whole space or a bounded domain.

We now recall the results in \cite{hw1} which yield the existence of
such a solution with the above properties precisely as
$\gamma>\f{N}{2}$, for $N=2,3$. We state the following theorem:

\begin{Theorem}[\bf{The periodic case}]\label{t1}
Assume that $\{(\r_\i, \u_\i, \H_\i)\}_{\i>0}$ is a sequence of weak
solutions to the compressible MHD equations \eqref{c} in the
periodic domain $\mathbb{T}$ with  initial data $\{(\r_\i^0,
\u_\i^0, \H_\i^0)\}_{\i>0}$,  satisfying the conditions
\eqref{6}-\eqref{8} and $\gamma>\f{N}{2}$, $N=2,3$. Also assume that
$(\u, \H)\in [L^2([0,T]; H^1(\mathbb{T}))\cap L^\infty([0,T];
L^2(\mathbb{T}))]^2$ is a weak solution to the incompressible MHD
equations \eqref{in} with initial data $\u|_{t=0}=P\u_0$ and $\H|_{t=0}=\H_0$. Then, for any finite number $T$, up to a subsequence, the global  weak
solutions $\{(\r_\i, \u_\i, \H_\i)\}_{\i>0}$ converge to $(\u,
\H)$. More precisely, as $\i\to 0$,
$$\r_\i \textrm{ converges to 1 in } C([0,T]; L^\gamma(\O));$$
$$P\u_\i \textrm{ converges strongly to } \u \textrm{ in } L^2([0,T]; L^p(\mathbb{T})), \textrm{ for all } 1\le p<\f{2N}{N-2};$$
$$Q\u_\i \textrm{ converges weakly to } 0 \textrm{ in } L^2([0,T]; H^1(\mathbb{T}));$$
$$\H_\i \textrm{ converges to } \H \textrm{ strongly in } L^2([0,T]; L^2(\mathbb{T})) \textrm{ and weakly in } L^2([0,T]; H^1(\mathbb{T})),$$
where, for convenience we will denote $\infty$ by $\f{2N}{N-2}$ if
$N=2$ in this paper.
\end{Theorem}

\medskip

\subsection{The whole space case}

Next, we turn to the whole space case. For the convenience of
presentation, we only discuss the case when $a=1$. In order to
define weak solutions in the whole space, the following special type
of Orlicz spaces $L^p_q(\O)$ are needed (see Appendix A in
\cite{L2}):
$$L^p_q(\O)=\left\{f\in L^1_{loc}(\O): \; f\chi_{\{|f|<\eta\}}\in L^q(\O), \;
 f\chi_{\{|f|\geq\eta\}}\in L^p(\O), \textrm{ for some } \eta>0\right\},$$
 where $\chi$ denotes the characteristic function of a set.
 We consider a sequence of
weak solutions $\{(\r_\i, \u_\i, \H_\i)\}_{\i>0}$ in the whole space
$\R^N$ with initial data $\{(\r_\i^0, \u_\i^0, \H_\i^0)\}_{\i>0}$,
satisfying the same conditions \eqref{6} and
\eqref{8} as in the periodic case. In addition, the weak solutions $\{(\r_\i,
\u_\i, \H_\i)\}_{\i>0}$ satisfy the following conditions at
infinity:
$$\r_\i\rightarrow 1,\qquad \u_\i\rightarrow 0,
\qquad\H_\i\rightarrow 0, \quad \textrm{ as } |x|\rightarrow \infty,$$
and
\begin{equation}\label{25}
\f{1}{2}\int_{\R^N}\left(\r_{\i}^0|\u_{\i}^0|^2+|\H_{\i}^0|^2\right)dx
+\f{a}{\i^2(\gamma-1)}\int_{\R^N}\left((\r_{\i}^0)^{\gamma}
-\gamma\r_{\i}^0+(\gamma-1)\right)dx\le C.
\end{equation}
As pointed out in \cite{hw1}, one can show that for any fixed
$\i>0$, there exists a global weak solution $(\r_\i, \u_\i, \H_\i)$
to the compressible MHD equations \eqref{c} defined by
$$\r_\i-1\in L^\infty([0,T]; L^\gamma_2(\R^N)),$$
$$\sqrt{\r_\i}\u_\i\in L^\infty([0,T]; L^2(\R^N)),$$
$$\nabla\u_\i\in L^2([0,T]; L^2(\R^N)),$$
$$\H_\i\in L^2([0,T]; H^1(\R^N))\cap L^\infty([0,T]; L^2(\R^N)),$$
satisfying, in addition,
$$\r_\i\u_\i\in C([0,T]; L^{2\gamma/(\gamma+1)}_{loc}(\R^N)),$$
$$\r_\i\in C([0,T]; L^p_{loc}(\R^N)),$$
if $1\le p<\gamma$ for all finite number $T$.

Now we are ready to state our result in the whole space as follows.

\begin{Theorem}[\bf{The whole space case}]\label{t2}
Assume that $\{(\r_\i, \u_\i, \H_\i)\}_{\i>0}$ is a sequence of weak
solutions to the compressible MHD equations \eqref{c} in the whole
space $\R^N$ with the initial data $\{(\r_\i^0, \u_\i^0,
\H_\i^0)\}_{\i>0}$, satisfying the conditions \eqref{6}, \eqref{8},
\eqref{25} and $\gamma>\f{N}{2}$, $N=2,3$. Also assume that $(\u,
\H)\in [L^2([0,T]; H^1(\R^N))\cap L^\infty([0,T]; L^2(\R^N))]^2$ is
a weak solution to the incompressible MHD equations \eqref{in} with
initial data $\u|_{t=0}=P\u_0$ and $\H|_{t=0}=\H_0$. Then, for any finite number $T$, up to a subsequence, the global  weak solutions $\{(\r_\i,
\u_\i, \H_\i)\}_{\i>0}$ converge to $(\u, \H)$. More precisely, as $\i\to 0$,
$$\r_\i \textrm{ converges to 1 in } C([0,T]; L^\gamma(\O));$$
$$P\u_\i \textrm{ converges strongly to } \u \textrm{ in } L^2([0,T]; L^p_{loc}(\R^N)), \textrm{ for all } 1\le p<\f{2N}{N-2};$$
$$Q\u_\i \textrm{ converges strongly to } 0 \textrm{ in } L^2([0,T]; L^q(\R^N)),\quad\textrm{for all}\quad 2<q<\f{2N}{N-2};$$
$$\H_\i \textrm{ converges to } \H \textrm{ strongly in } L^2([0,T]; L^2(\R^N)) \textrm{ and weakly in } L^2([0,T]; H^1(\R^N)).$$
\end{Theorem}

\medskip

\subsection{The bounded domain case}
The third case we will address in this paper is the incompressible
limit in a bounded domain $\O$. For the convenience of presentation, we
also only discuss the situation when $a=1$. In order to state
precisely our main theorem, we first introduce a geometrical
condition on $\O$ (cf. \cite{DGLM}). Let us consider the following
over-determined problem
\begin{equation}\label{30}
-\Delta\psi=\lambda\psi \quad\textrm{ in } \O,\qquad
\f{\partial\psi}{\partial\bf{n}}=0 \quad\textrm{ on }
\partial\O, \quad \textrm{ and } \psi \textrm{ is constant on }
\partial\O.
\end{equation}
A solution to \eqref{30} is said to be trivial if $\lambda=0$ and
$\psi$ is a constant. We   say that $\O$ satisfies the assumption
(A) if all the solutions to \eqref{30} are trivial. In the two
dimensional space, it is proved that every bounded, simply connected
open set $\O$ with Lipschitz boundary satisfies (A).

We consider a sequence of weak solutions $\{(\r_\i, \u_\i,
\H_\i)\}_{\i>0}$ in a bounded domain $\O$ with  initial data
$\{(\r_\i^0, \u_\i^0, \H_\i^0)\}_{\i>0}$ and boundary condition
\begin{equation} \label{bc1}
\u_\i|_{\partial\O}=0, \quad   \H_\i|_{\partial\O}=0,
\end{equation}
satisfying the same
conditions \eqref{6} and \eqref{8} as in the periodic case. And the
initial data of the weak solutions $\{(\r_\i, \u_\i, \H_\i)\}_{\i>0}$ satisfy
\begin{equation}\label{301}
\f{1}{2}\int_{\O}\left(\r_{\i}^0|\u_{\i}^0|^2+|\H_{\i}^0|^2\right)dx
+\f{a}{\i^2(\gamma-1)}\int_{\O}\left((\r_{\i}^0)^{\gamma}
-\gamma\r_{\i}^0+(\gamma-1)\right)dx \le C.
\end{equation}
As shown in \cite{hw1}, for any fixed $\i>0$, there exists a global
weak solution $(\r_\i, \u_\i, \H_\i)$ to the compressible MHD equations
\eqref{c} defined by
$$\r_\i\in L^\infty([0,T]; L^\gamma(\O)),$$
$$\sqrt{\r_\i}\u_\i\in L^\infty([0,T]; L^2(\O)),$$
$$\nabla\u_\i\in L^2([0,T]; H^1(\O)),$$
$$\H_\i\in L^2([0,T]; H^1(\O))\cap L^\infty([0,T]; L^2(\O)),$$
satisfying, in addition,
$$\r_\i\u_\i\in C([0,T]; L^{2\gamma/(\gamma+1)}(\O)),$$
$$\r_\i\in C([0,T]; L^p_{loc}(\O)),$$
if $1\le p<\gamma$ for all finite number $T$.

Our main result in bounded domains reads as follows.

\begin{Theorem}[\bf{The bounded domain case}]\label{t3}
Assume that $\{(\r_\i, \u_\i, \H_\i)\}_{\i>0}$ is a sequence of
weak solutions to the compressible MHD equations \eqref{c} in a
bounded domain $\O$ with  initial data $\{(\r_\i^0, \u_\i^0,
\H_\i^0)\}_{\i>0}$ and boundary condition \eqref{bc1},
satisfying the conditions \eqref{6},
\eqref{8}, \eqref{301} and $\gamma>\f{N}{2}$, $N=2,3$. Also assume
that $(\u, \H)\in [L^2([0,T]; H^1(\O))\cap L^\infty([0,T];
L^2(\O))]^2$ is a weak solution to the incompressible MHD
equations \eqref{in} with initial data $\u|_{t=0}=P\u_0$ and
$\H|_{t=0}=\H_0$ and boundary conditions $\u|_{\partial\O}=0$ and
$\H|_{\partial\O}=0$. Then for any finite number $T$, as $\i$ goes
to $0$, the global weak solutions $\{(\r_\i, \u_\i,
\H_\i)\}_{\i>0}$ converges to $(\u, \H)$. More precisely, as
$\i\to 0$,
$$\r_\i \textrm{ converges to 1 in } C([0,T]; L^\gamma(\O));$$
$$\u_\i \textrm{ converges to } \u \textrm{ weakly in } L^2(\O\times(0,T)) \textrm{ and strongly if } \O \textrm{ satisfies
(A)};$$
$$\H_\i \textrm{ converges to } \H \textrm{ strongly in } L^2([0,T]; L^2(\O)) \textrm{ and weakly in } L^2([0,T]; H^1(\O)).$$
\end{Theorem}

\begin{Remark}
In fact, we will split the eigenvectors
$\{\Psi_{k,0}\}_{k\in\mathbb{N}}$ of the Laplace equation with
Neumann boundary condition into two classes: those which are not
constant on $\partial\O$ will generate boundary layer and will be
quickly damped, thus converge strongly to 0; those which are
constant on $\partial\O$, for which no boundary layer forms, will
remain oscillating forever, and lead to only weak convergence.
Hence, if (A) is not satisfied, $\u_\i$ will in general only
converge weakly and not strongly to $\u$. In particular, in the
bounded, simply connected open set $\O\subset\R^2$ with Lipschitz
boundary, the boundary layer will always be generated, and hence $\u_\i$
will strongly converge to zero.
\end{Remark}

\bigskip

\section{The Periodic Case}

In this section, we will prove Theorem \ref{t1}.

\subsection{\textit{A priori} bounds and consequences}

We first deduce from \eqref{8} and from the conservation of mass
that we have for almost all $t\ge 0$,
\begin{equation}\label{101}
\begin{split}
&\f{1}{2}\int_{\mathbb{T}}\left(\r_{\i}|\u_{\i}|^2+|\H_{\i}|^2+\f{a}{\i^2(\gamma-1)}\left(\r_{\i}^\gamma
-\gamma\r_{\i}(\ov{\r_{\i}})^{\gamma-1}+(\gamma-1)(\ov{\r_{\i}})^\gamma\right)
\right)dx\\
&\qquad +\int_0^t\int_{\mathbb{T}}\left(\mu_{\i}|D\u_{\i}|^2+\lambda_{\i}(\Dv\u_{\i})^2+\nu_{\i}|\nabla\times\H_{\i}|^2\right) dxds\\
&\le \f{1}{2}\int_{\mathbb{T}}\left(\r^0_{\i}|\u^0_{\i}|^2+|\H^0_{\i}|^2+\f{a}{\i^2(\gamma-1)}\left((\r^0_{\i})^\gamma -\gamma\r^0_{\i}
(\ov{\r_{\i}})^{\gamma-1}+(\gamma-1)(\ov{\r_{\i}})^\gamma\right)\right)dx
\le C.
\end{split}
\end{equation}
From this inequality we see that $\r_{\i}|\u_{\i}|^2$,
$|\H_\i|^2$ and
$\f{1}{\i^2}\left(\r_\i^\gamma-\gamma\r_\i(\ov{\r_\i})^{\gamma-1}+(\gamma-1)(\ov{\r_\i})^\gamma\right)$
are bounded in $L^{\infty}([0,T]; L^1(\mathbb{T}))$ and that
$D\u_\i$ and $\nabla\times\H_\i$ are bounded in $L^2([0,T];
L^2(\mathbb{T}))$. In particular, we see that $\r_\i$ is bounded
in $L^\infty([0,T]; L^\gamma(\mathbb{T}))$ for all $T\in (0,\infty)$
due to the fact that for $\i$ small enough,
$\ov{\r_\i}\in(\f{1}{2},\f{3}{2})$ and thus for all $\dl>0$, there
exists some $\eta>0$ such that
\begin{equation}\label{100}
x^\gamma+(\gamma-1)(\ov{\r_\i})^\gamma-\gamma
x(\ov{\r_\i})^{\gamma-1}\ge \eta|x-\ov{\r_\i}|^\gamma \quad
\textrm{if}\quad |x-\ov{\r_\i}|\ge \dl,\quad x\ge 0.
\end{equation}
As in \cite{LM}  $\u_\i$ is bounded in
$L^2([0,T]; H^1(\mathbb{T}))$ for all $T\in (0,\infty)$. In fact, we
deduce from H\"{o}lder and Poinc\'{a}re's inequalities that we have
for all $T\in(0,\infty)$
\begin{equation*}
\int_0^T\int_{\mathbb{T}}\r_\i
\left|\u_\i-(2\pi)^{-N}\int_{\mathbb{T}}\u_\i dx\right|^2 dxdt
\le C\|\r_\i\|_{L^\infty([0,T]; L^\gamma)}\|D\u_\i\|^2_{L^2([0,T];
L^2)}\le C,
\end{equation*}
hence, in view of the above bound on $\r_\i|\u_\i|^2$, we get
\begin{equation*}
C\ge \int_0^T\int_{\mathbb{T}}\r_\i\left(\int_{\mathbb{T}}\u_\i dx\right)^2 dxdt
=\left(\int_{\mathbb{T}}\r_\i^0 dx\right)\int_0^T\left(
\int_{\mathbb{T}}\u_\i dx\right)^2dt.
\end{equation*}
Since \eqref{7} implies that $\r_\i^0$ converges to 1 in measure,
and hence, up to a subsequence, in $L^1(\mathbb{T})$, thus, we can
deduce a bound on $\u_\i$ in $L^2([0,T]; L^2)$ by using Poinc\'{a}re
inequality again. Indeed, we have
\begin{equation*}
\begin{split}
\int_0^T\int_{\mathbb{T}}|\u_\i|^2dxdt&\le
2\int_0^T\int_{\mathbb{T}}\left|\u_\i-(2\pi)^{-N}\int_{\mathbb{T}}\u_\i dx\right|^2dxdt+2(2\pi)^{-N}\int_0^T\left|\int_{\mathbb{T}}
\u_\i\right|^2dxdt\\&\le
C\left(1+\int_0^T\int_{\mathbb{T}}|\nabla\u_\i|^2dxdt\right)\le C.
\end{split}
\end{equation*}
From now on, we assume that, up to a subsequence,  $\u_\i$
converges weakly to some $\u$ in $L^2([0,T]; H^1(\mathbb{T}))$ for
all $T>0$. On the other hand, the bound on $\H_\i$ in
$L^\infty([0,T]; L^2(\mathbb{T}))$ and the bound on $\nabla\H_\i$
in $L^2([0,T]; L^2(\mathbb{T}))$, combining the following
Gagliardo-Nirenberg inequality $$\|\u\|_{L^{\f{8}{3}}([0,T];
L^4(\mathbb{T}))}\le \|\u\|^{\f{1}{4}}_{L^\infty([0,T];
L^2(\mathbb{T}))}\|\nabla\u\|_{L^2([0,T];
L^2(\mathbb{T}))}^{\f{3}{4}},$$ imply that $\H_\i$ is bounded in
$L^{8/3}([0,T]; L^4(\mathbb{T}))$, and also we can assume that
$\H_\i$ converges weakly to some $\H$ in $L^2([0,T];
H^1(\mathbb{T}))$ with $\Dv\H=0$. Finally, from the induction
equation in \eqref{c}, we see that
$$\partial_t\H_\i=\nabla\times(\u_\i\times\H_\i)-\nabla\times(\nu_\i\nabla\times\H_i)$$ is bounded in $L^{8/7}([0,T];
H^{-1}(\mathbb{T}))$, because the fact that $\u_\i$ is bounded in
$L^2([0,T]; L^4(\mathbb{T}))$ implies that $\u_\i\times\H_\i$ and
$\nu_\i\nabla\times\H_i$ are bounded in $L^{8/7}([0,T];
L^{2}(\mathbb{T}))$. Then the Aubin-Lions compactness Lemma (see
\cite{L2}) implies that $$\H_\i \rightarrow \H,\quad
\textrm{strongly in}\quad L^{8/7}([0,T];L^2(\mathbb{T})).$$
Moreover, this, combing with the uniform bound on $\H_\i$ in
$L^\infty([0,T]; L^2(\mathbb{T}))$, implies that $\H_\i$ converges
strongly to $\H$ in $L^{2}([0,T];L^2(\mathbb{T}))$. Therefore, by
a standard argument, we deduce that the limits $\u$ and $\H$
satisfy the induction equation in the sense of distributions, and
also the nonlinear term $(\nabla\times\H_\i)\times\H_\i$ in the
second equation of \eqref{c} converges to
$(\nabla\times\H)\times\H$ in the sense of distributions.

Next, we claim that $\r_\i$ converges to 1 in $C([0,T];
L^\gamma(\mathbb{T}))$. Indeed, in view of \eqref{101} and
\eqref{100}, we have
\begin{equation*}
\begin{split}
\sup_{t\ge 0}\int_{\mathbb{T}}|\r_\i-1|^\gamma dx
&\le
\dl^\gamma(2\pi)^N+C\sup_{t\ge 0}\left(\int_{\mathbb{T}}
\chi_{\{|\r_\i-1|\ge \dl\}}|\r_\i-\ov{\r_\i}|^\gamma dx\right)
+C|\ov{\r_\i}-1|^\gamma\\
&\le (2\pi)^N\dl^\gamma+\f{C\i^2}{\eta}+C|\ov{\r_\i}-1|^\gamma,
\end{split}
\end{equation*}
and we conclude the claim upon letting first $\i$ go to 0 and then $\dl$ go to
0.

Now, we show from the previous bounds that $\Dv\u_\i$ converges
weakly to 0 in $L^2([0,T]; L^2(\mathbb{T}))$ and that $P\u_\i$
converges to $\u=P\u$ strongly in $L^2([0,T]; L^2(\mathbb{T}))$,
and thus by Sobolev imbedding in $L^2([0,T]; L^q)$ for all $2\le
q<\f{2N}{N-2}$. These facts imply that $Q\u_\i$ converges weakly
to $0$ in $L^2([0,T]; H^1(\mathbb{T}))$. Indeed, since $\r_\i$
converges to $1$ in $C((0,\infty); L^\gamma(\mathbb{T}))$ and
$\gamma>\f{N}{2}$, we deduce from \eqref{c} that $\Dv\u_\i$
converges weakly to 0 in $L^2([0,T]; L^2(\mathbb{T}))$. The second
part is proven by observing first that we project \eqref{c} onto
divergence-free vector-fields:
\begin{equation}\label{10}
\partial_t P(\r_\i\u_\i)+P[\Dv(\r_\i\u_\i\otimes\u_\i)]-\mu_\i
\Delta P\u_\i=P((\nabla\times\H_\i)\times\H_\i).
\end{equation}
Noticing the fact that the operator $P$ is bounded in all Sobolev
space $W^{s,p}$ for all $s\in [0,\infty)$ and $1<p<\infty$ and the
preceding bounds, \eqref{10} yields a bound on $\partial_t
P(\r_\i\u_\i)$ in $L^1([0,T]; H^{-1}(\mathbb{T}))+L^2([0,T];
L^{1}(\mathbb{T}))+L^2([0,T]; H^{-1}(\mathbb{T}))$, hence, in
$L^1([0,T]; H^{-1}(\mathbb{T}))$. In addition, $P(\r_\i\u_\i)$ is
bounded in $L^\infty([0,T];
L^{\f{2\gamma}{\gamma+1}}(\mathbb{T}))\cap L^2([0,T];
L^r(\mathbb{T}))$ with
$$\f{1}{r}=\f{1}{\gamma}+\f{N-2}{2N}.$$

Next, we will need the following compactness Lemma (cf. Lemma 5.1
in \cite{L2}):
\begin{Lemma}
Let $g_n$, $h_n$ converge weakly to $g$, $h$ respectively in
$L^{p_1}(0,T; L^{p_2})$, $L^{q_1}(0,T; L^{q_2})$ where $1\le p_1,
p_2\le\infty,$
$$\f{1}{p_1}+\f{1}{q_1}=\f{1}{p_2}+\f{1}{q_2}=1.$$
Assume in addition that
$$\f{\partial g_n}{\partial t}\quad\textrm{is bounded in}\quad
L^1(0,T; W^{-m,1})\quad\textrm{for some}\quad m\ge 0\quad
\textrm{independent of n},$$ and
$$\|h_n-h_n(\cdot+\xi, t)\|_{L^{q_1}(0,T; L^{q_2})}\rightarrow
0\quad\textrm{as}\quad |\xi|\rightarrow 0,\quad\textrm{uniformly
in n}.$$ Then $g_nh_n$ converges to $gh$ in the sense of
distributions in $\O\times (0,T)$.
\end{Lemma}

Applying this lemma with the previous bounds, we deduce that
$P(\r_\i\u_\i)\cdot P\u_\i$ converges in the sense of
distributions to $|\u|^2$. We then conclude easily that $P\u_\i$
converges in $L^2([0,T]; L^2(\mathbb{T}))$ to $\u$ upon using the
weak convergence of $P\u_\i$ to $\u$ in $L^2([0,T];
L^2(\mathbb{T}))$ and remarking that we have
\begin{equation*}
\left|\int_0^T\int_{\mathbb{T}}\left(|P\u_\i|^2-P(\r_\i\u_\i)\cdot
P\u_\i\right)dxdt\right|\le
C\|\r_\i-1\|_{C([0,T];L^\gamma)}\|\u_\i\|^2_{L^2([0,T]; L^s)},
\end{equation*}
with $s=\f{2\gamma}{\gamma-1}<\f{2N}{N-2}$ since $\gamma>\f{N}{2}$.

We conclude this first step by showing the following bounds valid
for all $R\in (1,\infty)$
\begin{equation}\label{11}
\begin{cases}
&\|\varphi_\i\|_{L^\infty([0,T]; L^2(\mathbb{T}))}\le C\quad
\mathrm{if}\quad
\gamma\ge 2,\\
&\|\varphi_\i\chi_{\{\r_\i<R\}}\|_{L^\infty([0,T]; L^2(\mathbb{T}))}\le
C\quad \mathrm{if}\quad
\gamma< 2,\\
&\|\varphi_\i\chi_{\{\r_\i\ge R\}}\|_{L^\infty([0,T];
L^\gamma(\mathbb{T}))}\le C\i^{\f{2}{\gamma}-1}\quad
\mathrm{if}\quad
\gamma< 2,\\
\end{cases}
\end{equation}
where we denote the density fluctuation by
$\varphi_\i=\f{1}{\i}(\r_\i-\ov{\r_\i})$. These bounds are deduced
immediately from the following straightforward inequalities:  for some $\nu>0$
and for all $x\ge 0$,
\begin{equation}\label{104}
\begin{cases}
&x^\gamma-1-\gamma(x-1)\ge \nu|x-1|^2\quad \mathrm{if}\quad\gamma\ge
2,\\
&x^\gamma-1-\gamma(x-1)\ge \nu|x-1|^2\quad \mathrm{if}\quad\gamma<
2\quad\mathrm{and}\quad x\le R,\\
&x^\gamma-1-\gamma(x-1)\ge \nu|x-1|^\gamma\quad
\mathrm{if}\quad\gamma< 2\quad\mathrm{and}\quad x\ge R.
\end{cases}
\end{equation}

\subsection{The weak convergence of $Q\u$}

The proof of the weak convergence of $Q\u$ is similar to that in Lions-Masmoudi \cite{LM}, thus we only briefly describe the main idea from \cite{LM}. We provide here first a formal proof of
the passage to the limit, next the main difficulty,  and finally the
strategy of proof used in order to circumvent that difficulty.

We thus begin by an informal proof. It is not difficult to check
that the main difficulty with the passage to the limit lies with the
term $\Dv(\r_\i\u_\i\otimes\u_\i)$ and more precisely with the term
$\Dv(\r_\i Q(\u_\i)\otimes Q\u_\i)$ since the strong convergence of
$P\u_\i$. Formally, this term should not create an obstruction since
in view of the continuity equation in \eqref{1}, we can rewrite the
term $\partial_t(\r_\i\u_\i)+\Dv(\r_\i\u_\i\otimes\u_\i)$ as
$\r_\i\partial_t\u_\i+\r_\i(\u_\i\cdot\nabla)\u_\i$, which
corresponds to the term $\partial_t\u+(\u\cdot\nabla)\u$ in
incompressible MHD equations \eqref{in}. Next, the dangerous term
$[(Q\u_\i)\cdot\nabla]Q\u_\i$ can be incorporated in the pressure
$p$ at the limit since $Q\u_\i=\nabla\psi_\i$ for some $\psi_\i$,
and then
$$[(Q\u_\i)\cdot\nabla]Q\u_\i=\nabla\left|\f{1}{2}\nabla\psi_\i\right|^2.$$

Next, we need to write down rigorously the proof of the convergence.
First, we introduce the following group $\{\mathcal{L}(t), t\in R\}$
defined by $e^{tL}$ where $L$ is the operator defined on
$\mathcal{D}'_0\times(\mathcal{D}')^N$, where
$\mathcal{D}'_0=\{\phi\in \mathcal{D}', \int\phi=0\}$, by:
\begin{equation}\label{12}
L\binom{\phi}{\upsilon}=-\binom{\Dv\upsilon}{b\nabla\phi},\quad\textrm{
for}\quad b>0.
\end{equation}
We remark that $e^{tL}$ is an isometry on each $H^s\times(H^s)^N$
for all $s\in R$ and for all $t$, endowed with the norm
$\|(\phi,\upsilon)\|=(\|\phi\|_{H^s}^2+\f{1}{b}\|\upsilon\|^2_{H^s})^{1/2}$.
For details, we refer the reader to \cite{LM}. For convenience, in
the sequel, we will denote by $\mathcal{L}_1$ ($\mathcal{L}_2$)
the first (the second) component of the operator $\mathcal{L}$,
respectively.

We next claim that
$\mathcal{L}(-\f{t}{\i})\binom{\varphi_\i}{Q(\r_\i\u_\i)}$ is
relatively compact in $L^2([0,T]; H^{-n})$ for some $n\in(0,1)$. To
this end, we need to prove first that
$\binom{\varphi_\i}{Q(\r_\i\u_\i)}$ is bounded in $L^2([0,T];
H^{-s})$ for some $s\in(0,1)$ and that
$\partial_t\{\mathcal{L}\left(-\f{t}{\i}\right)\binom{\varphi_\i}{Q(\r_\i\u_\i)}\}$
is bounded in $L^2([0,T]; H^{-r})$ for some $r>0$ large enough. Our
claim then follows from Aubin-Lions compactness lemma by choosing n
in $(s,1)$.

Since we know that $\varphi_\i$ is bounded in $L^\infty([0,T]; L^p)$
where $p=\min(2,\gamma)$, by Sobolev imbedding theorems, we know
that $\varphi_\i$ is bounded in $L^2([0,T];H^{-s})$ for some
$s\in(0,1]$. And, we also deduce from the previous subsection that
$\r_\i\u_\i$ and thus $Q(\r_\i\u_\i)$ is bounded in $L^2([0,T];
L^q)$ with
$$\f{1}{q}=\f{1}{\gamma}+\f{N-2}{2N}.$$
 Therefore,
$\mathcal{L}\left(-\f{t}{\i}\right)\binom{\varphi_\i}{Q(\r_\i\u_\i)}$
is bounded in $L^2([0,T]; H^{-s})$ for some $s\in(0,1)$.

In order to get the uniform bound on
$\partial_t\{\mathcal{L}\left(-\f{t}{\i}\right)\binom{\varphi_\i}{Q(\r_\i\u_\i)}\}$,
we project the second equation of \eqref{c} into the space of
gradient vector-fields and we find
\begin{equation}\label{102}
\begin{split}
&\partial_tQ(\r_\i\u_\i)+Q[\Dv(\r_\i\u_\i\otimes\u_\i)]-(\mu_\i+\lambda_\i)\nabla\Dv\u_\i+\\
&\qquad\f{a}{\i^2}\nabla\left(\r_\i^\gamma-\gamma\r_\i(\ov{\r_\i})^{\gamma-1}+(\gamma-1)(\ov{\r_\i})^\gamma\right)+\f{a\gamma(\ov{\r_\i})^{\gamma-1}}
{\i^2}\nabla(\r_\i-\ov{\r_\i})\\
&=Q[(\nabla\times\H_\i)\times\H_\i].
\end{split}
\end{equation}
Hence, we can write the first equation \eqref{c} and \eqref{102} as
\begin{equation*}
\begin{split}
&\i\f{\partial\varphi_\i}{\partial t}+\Dv Q(\r_\i\u_\i)=0,\\
&\i\f{\partial Q(\r_\i\u_\i)}{\partial t}+b\nabla\varphi_\i=\i F_\i,
\end{split}
\end{equation*}
where $b=a\gamma (\ov{\r_\i})^{\gamma-1}$, and
\begin{equation*}
\begin{split}
F_\i=&(\mu_\i+\lambda_\i)\nabla\Dv\u_\i-Q[\Dv(\r_\i\u_\i\otimes\u_\i)] \\
 &-a\nabla\left(\f{1}{\i^2}(\r_\i^\gamma-\gamma\r_\i(\ov(\r_\i))^{\gamma-1}
+(\gamma-1)(\ov{\r_\i})^\gamma)]+Q[(\nabla\times\H_\i)\times\H_\i\right).
\end{split}
\end{equation*}
Since $b$ goes to $a\gamma$ as $\i$ goes to 0. we will ignore the
dependence of $b$ on $\i$ hereafter.
Now, we set
$$\psi_\i(t)=\mathcal{L}_1(-\f{t}{\i})\binom{\varphi_\i}{Q(\r_\i\u_\i)},\quad
m_\i(t)=\mathcal{L}_2(-\f{t}{\i})\binom{\varphi_\i}{Q(\r_\i\u_\i)},$$
then we have
\begin{equation*}
\begin{split}
\f{\partial}{\partial
t}\binom{\psi_\i}{m_\i}&=\mathcal{L}\left(-\f{t}{\i}\right)\left\{\f{\partial}{\partial
t}\binom{\varphi_\i}{Q(\r_\i\u_\i)}+\f{1}{\i}\binom{\Dv
Q(\r_\i\u_\i)}{b\nabla\varphi_\i}\right\}\\
&=\mathcal{L}\left(-\f{t}{\i}\right)\binom{0}{F_\i},
\end{split}
\end{equation*}
where $F_\i$ is bounded in $L^2([0,T];
H^{-1}(\mathbb{T}))+L^2([0,T]; W^{-1-\dl,1}(\mathbb{T}))$ for all
$\dl>0$, and hence, is bounded in $L^2([0,T]; H^{-r}(\mathbb{T}))$
for all $r>\f{N}{2}+1$. Thus $\f{\partial}{\partial
t}\binom{\psi_\i}{m_\i}$ is bounded in $L^2([0,T];
H^{-r}(\mathbb{T}))$.

We deduce from the compactness of $(\psi_\i, m_\i)$ that we may
assume without loss of generality that $(\psi_\i, m_\i)$ converges
in $L^2([0,T]; H^{-n})$ to some $(\psi, m)$. Since $Pm_\i=0$, we
also have $Pm=0$. Similarly, $\int\psi=0$. Hence, we have
\begin{equation}
\binom{\varphi_\i}{Q(\r_\i\u_\i)}=\mathcal{L}\left(\f{t}{\i}\right)\binom{\varphi}{m}+r_\i,\quad
r_\i\rightarrow 0\quad\mathrm{in}\quad L^2([0,T];
H^{-n})\quad\mathrm{as}\quad\i\rightarrow 0.
\end{equation}

Finally, following the argument of Step 4 in Section 3 in
\cite{LM}, one can show that $\varphi$, $m\in L^2([0,T];
L^2(\mathbb{T}))$ and $\Dv(\r_\i\u_\i\otimes\u_\i)-\Dv(v_\i\otimes
v_\i)$ converges to $\Dv(\u\otimes\u)$ in the sense of
distributions, where
$v_\i=\mathcal{L}_2(\f{t}{\i})\binom{\varphi}{m}$. Moreover,
following the argument of Step 5 in Section 3 in \cite{LM}, we can
show that $\Dv(v_\i\otimes v_\i)$ converges to a distribution
which is a gradient. Note that the magnetic field does not affect
the argument of convergence of
$\Dv(\r_\i\u_\i\otimes\u_\i)-\Dv(v_\i\otimes v_\i)$ and
$\Dv(v_\i\otimes v_\i)$ because the magnetic field $\H$ does not
affect the integrability of $F_\i$ based on our estimates, thus we only state those
convergence results without proof. We refer the  reader
to \cite{LM} for details.

This finishes the proof of Theorem \ref{t1}.

\bigskip

\section{The Whole Space Case}

In this section, we prove Theorem \ref{t2}. The idea is taken
from \cite{DG}. Before we start, we introduce homogeneous Sobolev
spaces for $1<p<\infty$ and $s\in R$ defined as usual by
$$\dot{W}^{s,p}(\R^N)=(-\Delta)^{-s/2}L^p(\R^N)\quad\mathrm{and}\quad
\dot{H}^s(\R^N)=\dot{W}^{s,2}(\R^N),$$ where $\Delta$ is the Laplace
operator.

Let us denote by $\zeta\in C_0^\infty(\R^N)$  a smoothing kernel such
that $\zeta\ge 0$, $\int_{\R^N}\zeta dx=1$, and define
$\zeta_\alpha(x)=\alpha^{-N}\zeta(x/\alpha)$. The following
estimate will be useful in this section (cf. \cite{DG}):
\begin{equation}\label{13}
\|f-f\ast\zeta_\alpha\|_{L^q}\le C\alpha^{1-\sigma}\|\nabla
f\|_{L^2},\quad \textrm{for all}\quad f\in\dot{H}^1,
\end{equation}
where $$q\in\left[2,\f{2N}{N-2}\right)\quad
\textrm{and}\quad\sigma=N\left(\f{1}{2}-\f{1}{q}\right),$$ and for
$1<p_2<p_1<\infty$, $s\ge 0$ and $\alpha\in(0,1)$, we have
\begin{equation}\label{14}
\|g\ast\zeta_\alpha\|_{L^{p_1}(\R^N)}\le
C\alpha^{-s-N(1/p_2-1/p_1)}\|g\|_{W^{-s,p_2}(\R^N)}.
\end{equation}

\subsection{\textit{A priori} estimates and consequences}

Most of the arguments developed in the periodic case can be adapted to
the whole space case. First, we obtain bounds on $D\u_\i$ in
$L^2([0,T]; L^2(\R^N))$ , on $\nabla\times\H_\i$ in $L^2([0,T];
L^2(\R^N))$ and on $\r_\i|\u_\i|^2$, and
$\f{1}{\i^2}(\r_\i^\gamma+(\gamma-1)-\gamma\r_\i)$ in
$L^\infty([0,T]; L^1(\R^N))$. The bound on $\u_\i$ in $L^2([0,T];
L^2(\R^N))$ follows from \eqref{104} and the following observation:
\begin{equation}\label{15}
\int_{R^N}\left(\f{1}{\i^2}|\r_\i-1|^2\chi_{\{|\r_\i-1|\le
1/2\}}+\f{1}{\i^2}|\r_\i-1|^\gamma \chi_{\{|\r_\i-1|\ge 1/2\}}\right)\le C,
\end{equation}
and thus, in particular,
\begin{equation*}
\begin{split}
\int_{\R^N}|\u_\i|^2 dx&\le C+\int_{\R^N}|\u_\i|^2\chi_{\{\r_\i\le 1/2\}} dx\\
&\le C+\left(\int_{\R^N}\chi_{\{\r_\i\le 1/2\}}dx\right)^{1/\gamma}
\left(\int_{\R^N}|\u_\i|^{2\gamma'}dx\right)^{1/\gamma'}\\
&\le C\left(1+(\mathrm{meas}(|\r_\i-1|\ge
1/2))^{1/\gamma}\|\u_\i\|_{L^2}^{2\theta}\|D\u_\i\|_{L^2}^{2(1-\theta)}\right)\\
&\le
C\left(1+\i^{2/\gamma}\|\u_\i\|_{L^2}^{2\theta}\|D\u_\i\|_{L^2}^{2(1-\theta)}\right),
\end{split}
\end{equation*}
where $$\f{\theta}{2}+(1-\theta)\f{N-2}{2N}=\f{1}{2\gamma'}.$$
 We then complete  the proof of our claim using the bound on $D\u_\i$
in $L^2([0,T]; L^2(\R^N))$ and the classical Young's inequality.
Moreover, if we define the density fluctuation as
$$\varphi_\i=\f{\r_\i-1}{\i},$$
then, it is bounded uniformly in $\i$ in
$L^\infty([0,T];L^\kappa_2)$ with $\kappa=\min\{2,\gamma\}$.
Furthermore, if we write
$$\u_\i=\u_\i^1+\u_\i^2,$$
 where
$$\u_\i^1=\u_\i\chi_{\{|\r_\i-1|\le 1/2\}}\quad\textrm{and}\quad\u_\i^2=\u_\i\chi_{\{|\r_\i-1|>
1/2\}},$$ then, we have $$\sup_{t\ge 0}\int_{\R^N}|\u_\i^1|^2dx\le
2\sup_{t\ge 0}\int_{\R^N}\r_\i|\u_\i|^2dx\le C,$$ and for $p<\kappa$
when $N=2$, $p=2\kappa/3$ if $N=3$,
\begin{equation*}
\begin{split}
\int_{\R^N}|\u_\i^2|^2dx&\le
C\int_{\R^N}|\r_\i-1|^p\chi_{\{|\r_\i-1|>1/2\}}|\u_\i|^2dx\\
&\le C\|(\r_\i-1)\chi_{\{|\r_\i-1|>1/2\}}\|^p_{L^\infty([0,T];
L^\kappa(\R^N))}\|\u_\i\|^2_{L^{2\kappa/(\kappa-p)}}\\
&\le
C\i^{2p/\kappa}\|\u_\i\|^{2-pN/\kappa}_{L^2(\R^N)}\|\nabla\u_\i\|^{pN/\kappa}_{L^2(\R^N)},
\end{split}
\end{equation*}
hence, by Young's inequality, $\u_\i^1$ is bounded in
$L^\infty([0,T]; L^2(\R^N))$ and $\u_\i^2\i^{-\beta}$ is bounded in
$L^2([0,T]; L^2(\R^N))$, where $\beta\in (0,1)$ if $N=2$ and
$\beta=2/3$ if $N=3$.

Recalling that $\gamma>N/2$, we deduce that $\u_\i$ is bounded in
$$L^2([0,T]; L^4(\R^N)\cap L^{2\gamma/(\gamma-1)}(\R^N)).$$
Hence, we have
$$\|\varphi_\i\u_\i\|_{L^2([0,T];
L^{4/3}(\R^N)+L^{2\kappa/(\kappa+1)}(\R^N)}\le C.$$ Therefore, using
Sobolev's imbedding, we deduce
$$\|\varphi_\i\u_\i\|_{L^2([0,T]; H^{-1}(\R^N))}\le C.$$

Finally, we already know that $\varphi^0_\i$ is bounded in
$L^\kappa_2(\R^N)$, hence in $H^{-1}(\R^N)$, since $\gamma>N/2$. On
the other hand, $m_\i^0$ can be rewritten as
$$m_\i^0=\f{m_\i^0}{\sqrt{\r_\i^0}}\sqrt{\r_\i^0}\chi_{\{|\r_\i^0-1|\le
1/2\}}+\f{m_\i^0}{\sqrt{\r_\i^0}}\f{\sqrt{\r_\i^0}}{\sqrt{|\r_\i^0-1|}}\sqrt{|\r_\i^0-1|}\chi_{\{|\r_\i^0-1|>
1/2\}}.$$ This implies that $m_\i^0$ is bounded in
$L^2(\R^N)+L^{2\kappa/(\kappa+1)}(\R^N)$, and hence in
$H^{-1}(\R^N)$. Therefore, $\binom{\varphi_\i^0}{m_\i^0}$ is bounded
in $H^{-1}(\R^N)$ uniformly in $\i$.

\medskip

\subsection{Strong convergence of $Q\u_\i$ to $0$}

We  now  prove that the gradient part of the velocity $Q\u_\i$
converges strongly to 0. More precisely, we claim that $Q\u_\i$
converges strongly to 0 in $L^2([0,T]; L^p(\R^N))$ for all
$p\in(2,\f{2N}{N-2})$ . Indeed, let us first observe that the
compressible MHD equations can be rewritten in terms of the density
fluctuation $\varphi_{\i}$, the momentum
$m_\i=\r_\i\u_\i$ and $\phi_\i=\binom{\varphi_\i}{m_\i}$ as follows
$$\partial_t\phi_\i+\f{L\phi_\i}{\i}=F^1_\i+F^2_\i,$$
where the wave operator L is defined on $(\mathcal{D}'(\R^N))^{N+1}$
with values in $(\mathcal{D}'(\R^N))^{N+1}$ by
$$L\phi=\binom{\Dv \,m}{\nabla\psi},\quad \textrm{with} \quad
\phi=\binom{\psi}{m},$$ and
\begin{equation*}
\begin{split}
F^1_\i&=\binom{0}{\mu_\i\Delta\u^1_\i+\lambda_\i\nabla\Dv\u^1_\i-\Dv(m_\i\otimes\u_\i)-
\f{a}{\i^2}\nabla(\r_\i^\gamma-1-\gamma(\r_\i-1))+(\nabla\times\H_\i)\times\H_\i},\\
F^2_\i&=\binom{0}{\mu_\i\Delta\u^2_\i+\lambda_\i\nabla\Dv\u^2_\i}.
\end{split}
\end{equation*}

Using \textit{Duhamel}'s formula, we deduce that
$$Q\phi_\i(t)=\mathcal{L}\left(\f{t}{\i}\right)Q\phi^0_\i+\int_0^t\mathcal{L}\left(\f{t-s}{\i}\right)(QF^1_\i(s)+QF^2_\i(s))ds.$$
Here we used the fact that $Q$ and $\mathcal{L}$ commute, since $Q$
and $L$ do.

At this stage, the following Strichartz's estimates from \cite{DG}
are useful:

\begin{Lemma}\label{l1}
For all $s\ge 0$, we have
\begin{equation}\label{16}
\left\|\mathcal{L}\left(\f{t}{\i}\right)Q\psi_0\right\|_{L^q((0,\infty);
W^{-s-\sigma,p}(\R^N))}\le C\i^{1/q}\|\psi_0\|_{H^{-s}(\R^N)},
\end{equation}
\begin{equation}\label{17}
\left\|\int_0^t\mathcal{L}\left(\f{t-s}{\i}\right)Q\psi(s)ds\right\|_{L^q([0,T];
W^{-s-\sigma,p}(\R^N))}\le
C(1+T)\i^{1/q}\|\psi\|_{L^q([0,T];H^{-s}(\R^N))},
\end{equation}
for all $(p,q)\in(2,\infty)\times(2,\infty)$ and
$\sigma\in(0,\infty)$ such that
\begin{equation}\label{18}
\f{2}{q}=(N-1)\left(\f{1}{2}-\f{1}{p}\right)\quad\mathrm{and}\quad\sigma
q=\f{N+1}{N-1}.
\end{equation}
\end{Lemma}

\bigskip

Now, we choose $p\in(2,\f{2N}{N-2})$, $q\in(2,\infty)$ and
$\sigma\in(0,\infty)$ given by \eqref{18}. One can deduce that
$$|Q\u_\i|\le
|Q\u_\i-Q\u_\i\ast\zeta_\alpha|+\i|Q(\u_\i\varphi_\i)\ast\zeta_\alpha|+|Qm_\i\ast\zeta_\alpha|.$$
Hence,
\begin{equation*}
\begin{split}
\|Q\u_\i\|_{L^2([0,T]; L^p(\R^N))}&\le
C\alpha^{1-N(1/2-1/p)}\|\nabla\u_\i\|_{L^2([0,T]; L^2(\R^N))}\\
&\quad+\i \alpha^{-1-N(1/2-1/p)}\|\varphi_\i\u_\i\|_{L^2([0,T];
H^{-1}(\R^N))}\\&\quad+\|Qm_\i\ast\zeta_\alpha\|_{L^2([0,T];
L^p(\R^N))}.
\end{split}
\end{equation*}
From the estimates in previous subsection, we know that $F^1_\i$ is
bounded in $L^\infty([0,T]; H^{-s_0})$ for all $s_0>N/2+1$. On the
other hand, we deduce from the uniform bound on $\u^2_\i
\i^{-\beta}$ in $L^2([0,T]; L^2)$ that $\i^{-\beta}F_\i^2$ is
bounded in $L^2([0,T]; H^{-2}(\R^N))$. Then, using Lemma \ref{l1},
we obtain,  for all $\eta>0$ small enough,
\begin{equation*}
\begin{split}
&\|Qm_\i\ast\zeta_\alpha\|_{L^2([0,T]; L^p(\R^N))}\\
&\le C_T\alpha^{-1-\sigma}
\left\|\mathcal{L}\left(\f{t}{\i}\right)\psi^0_\i\right\|_{L^q([0,T];
W^{-1-\sigma, p}(\R^N))}\\
&\qquad+C_T
\alpha^{-N/2-1-\sigma-\eta}\left\|\int_0^Tds\mathcal{L}\left(\f{t-s}{\i}\right)QF^1_\i(s)\right\|_{L^q([0,T];W^{-\eta-N/2-1-\sigma,
p}(\R^N))}\\
&\qquad+C
\alpha^{-2-N(1/2-1/p)}\left\|\int_0^tds\mathcal{L}\left(\f{t-s}{\i}\right)QF^2_\i(s)\right\|_{L^2([0,T];
H^{-2})}\\
&\le
C_T\alpha^{-1-\sigma}\i^{1/q}\|\psi^0_\i\|_{H^{-1}}+C_T\alpha^{-N/2-1-\sigma-\eta}\i^{1/q}\|F^1_\i\|_{L^\infty([0,T];
H^{-\eta-N/2-1})}\\&\qquad+C\alpha^{-2-N(1/2-1/p)}\i^\beta\|\i^{-\beta}F^2_\i\|_{L^2([0,T];
H^{-2})}.
\end{split}
\end{equation*}
Next, fixing $\alpha>0$ and letting $\i$ go to zero, we obtain
$$\limsup_{\i\rightarrow 0}\|Q\u_\i\|_{L^2([0,T]; L^p(\R^N))}\le C\alpha^{1-N(1/2-1/p)},$$
where $C$ is independent of $\i$ and $\alpha$.
Noticing that $1-N(1/2-1/p)>0$, we finally get, by letting
$\alpha\rightarrow 0$,
$$\limsup_{\i\rightarrow 0}\|Q\u_\i\|_{L^2([0,T]; L^p(\R^N))}=0.$$
This implies that $Q\u_\i$ strongly converges to $0$ in $L^2([0,T];
L^p(\R^N))$ for all $2<p<\f{2N}{N-2}$.

\medskip

\subsection{Strong convergences of $P\u_\i$ and $\H_\i$}

In the previous section, we proved the strong convergence of the
gradient part of the velocity to 0. In order to complete the proof
of Theorem \ref{t2}, we are left to show the convergence of the
incompressible part the velocity, $P\u_\i$, the convergence of the
density, and the convergence of the magnetic field. This can be done
by using the classical compactness arguments  in \cite{L2, LM},
or equivalently by looking at the time-regularity properties of
$P\u_\i$, see \cite{DG}. Indeed, following the argument in the
periodic case steps by steps, we  obtain the strong convergence
of $\r_\i$ to 1 in $C([0,T]; L^\gamma_{loc}(\R^N))$ and the weak
convergence of $P\u_\i$ to $\u$ in $L^2([0,T]; H^1(\R^N))$.
Moreover, we also can show that $P\u_\i$ converges to $\u$ in
$L^2([0,T]; L^2(B_R))$ for all $R\in(0,\infty)$. Here, we denote by
$B_R$ the open ball centered at 0 of radius $R$.

Finally, similarly as in the periodic case, the bound on $\H_\i$
in $L^\infty([0,T]; L^2(\R^N))$ and the bound on $\nabla\H_\i$ in
$L^2([0,T]; L^2(\R^N))$, combining Sobolev's inequality and
interpolation theorem, we know that $\H_\i$ is bounded in
$L^{8/3}([0,T]; L^4(\R^N))$, and also we can assume that $\H_\i$
converges weakly to some $\H$ in $L^2([0,T]; H^1(\R^N))$ with
$\Dv\H=0$. Finally, from the induction equation in \eqref{c}, we
deduce that $\partial_t\H_\i$ is bounded in $L^{8/7}([0,T];
H^{-1}(\R^N))$, due to the fact that $\u_\i$ is bounded in
$L^2([0,T]; L^4(\R^N))$. This property, combining Aubin-Lions
compactness Lemma, implies that $\H_\i$ converges strongly to $\H$
in $L^{8/7}([0,T];L^2_{loc}(\R^N))$. Moreover, the uniform bound
on $\H_\i$ in $L^\infty([0,T]; L^2(\R^N))$ implies that $\H_\i$
converges strongly to $\H$ in $L^{2}([0,T];L^2_{loc}(\R^N))$.
Therefore, by a standard argument, we deduce that the limits $\u$
and $\H$ satisfy the induction equation in \eqref{1} in the sense
of distributions, and also the nonlinear term
$(\nabla\times\H_\i)\times\H_\i$ in the second equation of
\eqref{c} converges to $(\nabla\times\H)\times\H$ in the sense of
distributions. For a detailed statement of the above argument, we
refer it to the argument surrounding the convergence of the
magnetic field in section 3.1.

The proof of Theorem \ref{t2} is complete.

\bigskip

\section{The Bounded Domain Case}
In this section, we will prove Theorem \ref{t3} by the spectral
analysis of the semigroup generated by the dissipative wave
operator. Before we start, we introduce the eigenvalues
$\{\lambda^2_{k,0}\}_{k\in\mathbb{N}}$ ($\lambda_{k,0}>0$) and the
eigenvectors $\{\Psi_{k,0}\}_{k\in\mathbb{N}}$ in $L^2(\O)$ with
zero mean value of the Laplace operator satisfying homogeneous Neumann
boundary conditions:
$$-\Delta\Psi_{k,0}=\lambda_{k,0}^2\Psi_{k,0}\quad\textrm{ in } \O,\qquad \f{\partial\Psi_{k,0}}{\partial\bf{n}}=0\quad \textrm{on }\partial\O.$$
Notice that, by Gram-Schmidt orthogonalization method, it is
possible to assume that $\{\Psi_{k,0}\}_{k\in\mathbb{N}}$ is a
orthonormal basis of $L^2(\O)$ and that up to a slight modification,
if $\lambda_{k,0}=\lambda_{l,0}$ and $k\neq l$, then
$$\int_{\partial\O}\nabla\Psi_{k,0}\cdot\nabla\Psi_{l,0}ds=0.$$

Next, we recall from the previous section that we can deduce
similarly that
$$\sup_{t\ge 0}\|\r_\i-1\|_{L^\gamma(\O)}\le C\i^{\kappa/\gamma}\quad
\textrm{and}\quad \sup_{t\ge 0}\|\r_\i-1\|_{L^\kappa(\O)}\le C\i,$$
where $\kappa=\min\{2,\gamma\}$. And similarly to the whole space
case, we will split
$$\u_\i=\u^1_\i+\u^2_\i, \quad\text{with}\quad
\u^1_\i=\u_\i\chi_{|\r_\i-1|\le 1/2},\quad
\u_\i^2=\u_\i\chi_{|\r_\i-1|>1/2},$$
 which satisfy
$$\sup_{t\ge 0}\int_\O|\u_\i^1|^2dx\le
2\sup_{t\ge0}\int_\O\r_\i|\u_\i|^2dx\le C;$$ and
$$\|\u^2_\i\|^2_{L^2(\O)}\le 2\int_\O|\r_\i-1||\u_\i|^2dx\le
C\i\|\u_\i\|^2_{L^{2\kappa/(\kappa-1)}(\O)}\le
C\i|\nabla\u_\i|^2_{L^2(\O)}.$$ Therefore, $\u^1_\i$ is bounded in
$L^\infty([0,T]; L^2(\O))$ whereas $\u_\i^2\i^{-1/2}$ is bounded in
$L^2(\O\times(0,T))$, and hence, $\u_\i$ is bounded in
$L^2(\O\times(0,T))$. Also, in this section, we denote the density
fluctuation by
$$\varphi_\i=\f{\r_\i-1}{\i},$$
and the momentum by $\bf{m}_\i=\r_\i\u_\i$.

\medskip

\subsection{Strong convergence of $P\u_\i$ and $\H_\i$}
Following the argument in the periodic case step by step, up to
the extraction of a subsequence, we then obtain the strong
convergence of $\r_\i$ to 1 in $C([0,T]; L^\gamma(\O))$, the strong
convergence of $P\u_\i$ to $\u=P\u$ in $L^2(\O\times(0,T))$, and the
weak convergence of $Q\u_\i$ to $0$ in $L^2([0,T]; H^1(\O))$. Thus,
the continuity equation in \eqref{c} holds in the sense of
distributions.

Similarly to the periodic case, the bound on $\H_\i$ in
$L^\infty([0,T]; L^2(\O))$ and the bound on $\nabla\H_\i$ in
$L^2([0,T]; L^2(\O))$, combining Sobolev's inequality and
interpolation theorem, we know that $\H_\i$ is bounded in
$L^{8/3}([0,T]; L^4(\O))$, and also we can assume that $\H_\i$
converges weakly to some $\H$ in $L^2([0,T]; H^1(\O))$ with
$\Dv\H=0$. Also, from the induction equation in \eqref{c}, we deduce
that $\partial_t\H_\i$ is bounded in $L^{8/7}([0,T]; H^{-1}(\O))$,
due to the fact that $\u_\i$ is bounded in $L^2([0,T]; L^4(\O))$.
This property, combining Aubin-Lions compactness Lemma, implies that
$\H_\i$ converges strongly to $\H$ in $L^{8/7}([0,T];L^2(\O))$.
Moreover, the uniform bound on $\H_\i$ in $L^\infty([0,T]; L^2(\O))$
implies that $\H_\i$ converges strongly to $\H$ in
$L^{2}([0,T];L^2(\O))$. Therefore, by a standard argument, we deduce
that the limits $\u$ and $\H$ satisfy the induction equation in
\eqref{1} in the sense of distributions, and also the nonlinear term
$(\nabla\times\H_\i)\times\H_\i$ in the second equation of \eqref{c}
converges to $(\nabla\times\H)\times\H$ in the sense of
distributions. Therefore, in order to prove Theorem \ref{t3}, it
only remains to study the convergence of the gradient part of the
velocity $Q\u_\i$.

\medskip

\subsection{The convergence of $Q\u_\i$}
The argument for the convergence of $Q\u_\i$ in this subsection
follows the lines in \cite{DGLM}, except the argument for the
magnetic field. For the reader's convenience and the completeness
of the argument, we provide the details here. For this purpose,
first, we discuss the spectral problem associated with the viscous
wave operator $L_\i$ in terms of eigenvalues and eigenvectors of
the inviscid wave operator $L$, where the wave operator $L$ and
$L_\i$ are defined on $\mathcal{D}'(\O)\times\mathcal{D}'(\O)^N$
by
$$L\binom{\psi}{\bf{m}}=\binom{\Dv \bf{m}}{\nabla\Psi},$$
and
$$L_\i\binom{\Psi}{\bf{m}}=L\binom{\Psi}{\bf{m}}+\i\binom{0}{\mu_\i\Delta
\bf{m}+\lambda_\i\nabla\Dv \bf{m}}.$$ The eigenvalues and
eigenvectors of $L$ read as follows
$$\phi_{k,0}^\pm=\binom{\Psi_{k,0}}{\bf{m}_{k,0}=\pm\f{\nabla\Psi_{k,0}}{\mathrm{i}\lambda_{k,0}}},$$
$$L\phi_{k,0}^\pm=\pm {\mathrm{i}}\lambda_{k,0}\phi_{k,0}^\pm\quad\textrm{in
}\O,\qquad {\bf{m}}_{k,0}^\pm\cdot{\bf{n}}=0 \quad\textrm{on
}\partial\O.$$

In the following steps, the following information on the
approximating eigenvalues and eigenvectors for the operator $L_\i$
is crucial:

\begin{Lemma}\label{l302}
Let $\O$ be a $C^2$ bounded domain in $R^N$ and let $k\ge 1$, $M\ge
0$. Then, there exists approximate eigenvalues
${\mathrm{i}}\lambda_{k,\i,M}^\pm$ and eigenvectors
$\phi_{k,\i,M}^\pm=\binom{\Psi_{k,\i,M}^\pm}{{\bf{m}}_{k,\i,M}^\pm}$
of $L_\i$ such that
$$L_\i\phi_{k,\i,M}^\pm={\mathrm{i}}\lambda^\pm_{k,\i,M}\phi_{k,\i,M}^\pm+R_{k,\i,M}^\pm,$$
with
$${\mathrm{i}}\lambda^\pm_{k,\i,M}=\pm
{\mathrm{i}}\lambda_{k,0}+{\mathrm{i}}\lambda^\pm_{k,1}\sqrt{\i}+O(\i),\quad
\textrm{ where } \mathcal{R}e({\mathrm{i}}\lambda^\pm_{k,1})\le 0,$$
and for all $1\le p\le \infty$, we have
$$\|R^\pm_{k,\i,M}\|_{L^p(\O)}\le C_p(\sqrt{\i})^{M+1/p}\quad
\textrm{and} \quad\|\phi^\pm_{k,\i,M}-\phi^\pm_{k,0}\|_{L^p(\O)}\le
C_p(\sqrt{\i})^{1/p}.$$
\end{Lemma}
\begin{proof}
For the construction in detail, we refer the readers to \cite{DGLM}.
\end{proof}

\begin{Remark}
Due to the construction in \cite{DGLM}, indeed, we have
$${\mathrm{i}}\lambda_{k,1}^\pm=-\f{1\pm
{\mathrm{i}}}{2}\sqrt{\f{\mu_\i}{2\lambda_{k,0}^3}}\int_{\partial_\O}\left|\nabla\Psi_{k,0}\right|^2ds.$$
\end{Remark}

\begin{Remark}
We notice that the first order term ${\mathrm{i}}\lambda^\pm_{k,1}$
clearly yields an instantaneous damping of the acoustic waves, as
soon as $\mathcal{R}e({\mathrm{i}}\lambda^\pm_{k,1})<0$. For this
reason, we define $I\subset \mathbb{N}$ to be the set of eigenvectors
$\Psi_{k,0}$ of the Laplace operator such that
$\mathcal{R}e({\mathrm{i}}\lambda^\pm_{k,1})<0$ and $J={\mathbb{N}}-I$.
Observe that when $k\in J$, we have $\lambda^\pm_{k,1}=0$. For those
indices, ${\bf{m}}_{k,0}^\pm$ identically vanishes on $\partial\O$
and therefore satisfies not only ${\bf{m}}_{k,0}^\pm\cdot\bf{n}=0$
but also ${\bf{m}}_{k,0}^\pm=0$ on $\partial\O$, hence no
significant boundary layer is created, and there is no enhanced
dissipation of energy in these layers.
\end{Remark}

Now, we can express $Q\u_\i$ in the terms of the orthonormal basis
$\left\{\f{\nabla\Psi_{k,0}}{\lambda_{k,0}}\right\}_{k\in\bf{N}}$ of
$L^2(\O)$ as
$$Q\u_\i=\sum_{k\in\bf{N}}\left(Q\u_\i,
\f{\nabla\Psi_{k,0}}{\lambda_{k,0}}\right)\f{\nabla\Psi_{k,0}}{\lambda_{k,0}},$$
where the notation $(\cdot, \cdot)$ stands for
$$(f(x),g(x))=\int_\O f(x)\ov{g(x)}dx.$$
We can split $Q\u_\i$ into two parts $Q_1\u_\i$ and $Q_2\u_\i$,
defined by
$$Q_1\u_\i=\sum_{k\in I}\left(Q\u_\i,
\f{\nabla\Psi_{k,0}}{\lambda_{k,0}}\right)\f{\nabla\Psi_{k,0}}{\lambda_{k,0}},\quad\textrm{and}\quad
Q_2\u_\i=\sum_{k\in J}\left(Q\u_\i,
\f{\nabla\Psi_{k,0}}{\lambda_{k,0}}\right)\f{\nabla\Psi_{k,0}}{\lambda_{k,0}},$$
which respectively correspond to damped terms and nondamped terms.
We will prove on one hand that $Q_1\u_\i$ converges strongly to $0$
in $L^2(\O\times(0,T))$, and on the other hand that $\mathrm{curl}
\Dv(Q_2 m_\i\otimes Q_2\u_\i)$ converges to 0 in the sense of
distributions, if $J\neq\emptyset$, which is equivalent to say that
$\Dv(Q_2 m_\i\otimes Q_2\u_\i)$ converges to a gradient in the sense
of distributions.

Let us observe that in view of the bound on $\u_\i$ in $L^2([0,T];
H^1_0(\O))$, the problem reduces to a finite number of terms.
Indeed, we have
$$\sum_{k>M}\int_0^T\left|\left(Q_i\u_\i,
\f{\nabla\Psi_{k,0}}{\lambda_{k,0}}\right)\right|^2dt\le
\f{C}{\lambda^2_{M+1}}|\nabla\u_\i|^2_{L^2(\O\times(0,T))}, \quad
i=1 \textrm{ or } 2.$$ Hence, recalling that
$\lambda_M\rightarrow\infty$ as $M\rightarrow\infty$, we only have
to prove that $(Q_1\u_\i, {\bf{m}}_{k,0}^\pm)$ converges strongly to
$0$ in $L^2(0,T)$ for any fixed $k$, and study the interaction of a
finite number of terms in $\Dv(Q_2\u_\i\otimes Q_2\u_\i)$. On the
other hand, we notice that
$$Q\u_\i=Q{\bf{m}}_\i-\i Q(\varphi_\i\u_\i),$$
 and
\begin{equation*}
\begin{split}
\i|(Q(\varphi_\i\u_\i), \nabla
\Psi_{k,0})|&=\i\left|\int_\O\varphi_\i\u_\i\cdot\nabla\Psi_{k,0}
dx\right|\\&\le\i\|\varphi_\i\|_{L^{\gamma}(\O)}\|\u_\i\|_{L^{\f{\gamma}{\gamma-1}}(\O)}\|\nabla\Psi_{k,0}\|_{L^\infty(\O)},
\end{split}
\end{equation*}
which goes to $0$ in $L^2(0,T)$ since $\gamma>N/2$. Hence, we are
led to study $(Q{\bf{m}}_\i, {\bf{m}}_{k,0}^\pm)$.

Denote $$\beta^\pm_{k,\i}=(\phi_\i(t), \phi^\pm_{k,0}),\quad\text{with}\quad
\phi_\i(t)=\binom{\varphi_\i}{\bf{m}_\i},$$ we observe that:
$$2(Q{\bf{m}}_\i,
{\bf{m}}_{k,0}^\pm)=\beta^\pm_{k,\i}-\beta^\mp_{k,\i},$$ so that it
suffices to consider the convergence properties of
$\beta^\pm_{k,\i}$ in $L^2(0,T)$. Also, we see from Lemma
\ref{l302} with $M=2$ that:
$$|(\phi_\i(t), \phi^\pm_{k,0}-\phi^\pm_{k,\i,2})|\le
C\i^{\alpha/2}\left(\|\varphi_\i\|_{L^\infty([0,T];
L^\kappa(\O))}+\|{\bf{m}}_\i\|_{L^\infty([0,T];
L^{\f{2\gamma}{\gamma+1}}(\O))}\right),$$ where
$$\alpha=\min\left\{1-\frac1\kappa, \frac12-\frac1{2\gamma}\right\},$$
hence we only have to
prove that $b_{k,\i}^\pm(t)=(\phi_\i(t), \phi^\pm_{k,\i,2})$
converges strongly to $0$ in $L^2(0,T)$ when $k\in I$, and study its
oscillations when $k\in J$.

Notice that $\phi_\i(t)=(\varphi_\i, {\bf{m}}_\i)$ solves
\begin{equation}\label{304}
\partial_t\phi_\i-\f{L_\i^*\phi_\i}{\i}=\binom{0}{g_\i},
\end{equation}
where $L_\i^*$ denotes the adjoint of $L_\i$ with respect to
$(\cdot,\cdot)$, and
$$g_\i=-\Dv({\bf{m}}_\i\otimes\u_\i)-\nabla \left[\f{(\r_{\i})^{\gamma}
-\gamma\r_{\i}+(\gamma-1)}{\i^2}\right]+(\nabla\times\H_\i)\times\H_\i.$$

Taking the scalar product of \eqref{304} with $\phi^\pm_{k,\i,2}$,
we obtain
\begin{equation}\label{305}
\f{d}{dt}b_{k,\i}^\pm(t)-\f{\ov{{\mathrm{i}}\lambda_{k,\i,2}^\pm}}{\i}b_{k,\i}^\pm(t)=c_{k,\i}^\pm(t),
\end{equation}
where $c_{k,\i}^\pm(t)=(g_\i,
{\bf{m}}^\pm_{k,\i,2})+\i^{-1}(\phi_\i, R^\pm_{k,\i,2})$.

\subsubsection{\bf{The case $k\in I$}}
From \eqref{305}, by Duhamel's principle, we deduce that
\begin{equation}\label{306}
b^\pm_{k,\i}(t)=b^\pm_{k,\i}(0)\exp^{\ov{{\mathrm{i}}\lambda^\pm_{k,\i,2}}t/\i}+\int_0^t
c_{k,\i}^\pm(s)\exp^{\ov{{\mathrm{i}}\lambda^\pm_{k,\i,2}}(t-s)/\i}ds.
\end{equation}
The first term in \eqref{306} is estimated as follows
$$\left\|b^\pm_{k,\i}(0)\exp^{\ov{{\mathrm{i}}\lambda^\pm_{k,\i,2}}t/\i}\right\|_{L^2(0,T)}\le
C\left\|b^\pm_{k,\i}(0)\exp^{\mathcal{R}e(\ov{{\mathrm{i}}\lambda^\pm_{k,1}})t/\sqrt{\i}}\right\|_{L^2(0,T)}\le
C\i^{1/4}.$$ In order to estimate the remaining term in \eqref{306},
we will use the following estimate: for any $1\le p,q\le\infty$ with
$\f{1}{q}+\f{1}{p}=1$, we have
\begin{equation}\label{307}
\left|\int_0^t\exp^{\ov{{\mathrm{i}}\lambda^\pm_{k,\i,2}}(t-s)/\i}a(s)ds\right|\le
\int_0^t\exp^{\mathcal{R}e(\ov{{\mathrm{i}}\lambda^\pm_{k,1}})(t-s)/\sqrt{\i}}|a(s)|ds\le
C\|a\|_{L^q(0,T)}{\i}^{\frac1{2p}}.
\end{equation}
We now write $|c_{k,\i}^\pm|\le c_1+c_2+c_3+c_4$, where
\begin{equation*}
\begin{split}
c_1(t)&=\left|\int_\O({\bf{m}}_\i\otimes\u_\i)(t)\cdot\nabla
{\bf{m}}^\pm_{k,\i,2}dx\right|,\\
c_2(t)&=\left|\int_\O \left[\f{(\r_{\i})^{\gamma}
-\gamma\r_{\i}+(\gamma-1)}{\i^2}\right](t)\Dv {\bf{m}}^\pm_{k,\i,2}dx\right|,\\
c_3(t)&=\i^{-1}|(\phi_\i, R^\pm_{k,\i,2})|,\\
c_4(t)&=\left|\int_\O(\nabla\times\H_\i)\times\H_\i\cdot\nabla
{\bf{m}}^\pm_{k,\i,2}dx\right|.
\end{split}
\end{equation*}
Observing that ${\bf{m}}_\i=\i\varphi_\i\u_\i+\u_\i$, we have:
\begin{equation*}
\begin{split}
c_1(t)&\le
\|{\bf{m}}_{k,\i,2}^\pm\|_{L^\infty(\O)}\|\u_\i^1+\u_\i^2\|_{L^2(\O)}\|\nabla\u_\i\|_{L^2(\O)}\\
&\quad+\i\|\varphi_\i\|_{L^\infty([0,T];
L^\kappa(\O))}\|(u_\i)^2\|_{L^{\kappa/(\kappa-1)}(\O)}\|\nabla
{\bf{m}}^\pm_{k,\i,2}\|_{L^\infty(\O)}\\
&\le C\|\u^1_\i\|_{L^\infty([0,T];
L^2(\O))}\|\nabla\u_\i\|_{L^2(\O)}+C\i^{1/2}\|\nabla\u_\i\|^2_{L^2(\O)}+C\i^{1/2}\|\nabla\u_\i\|_{L^2(\O)}.
\end{split}
\end{equation*}
The second term $c_2$ is estimated as
$$c_2(t)\le
C\|p_\i\|_{L^\infty([0,T];L^1(\O))}(\|\Psi_{k,\i,2}^\pm\|_{L^\infty(\O)}+\|R^\pm_{k,\i,2}\|_{L^\infty(\O)})\le
C.$$ Also, we estimate $c_3$ by
$$c_3(t)\le\f{1}{\i}\|R^\pm_{k,\i,2}\|_{L^{\kappa/(\kappa-1)}(\O)}\|\phi_\i\|_{L^\infty([0,T];
L^\kappa(\O))}\le C\i^{1/2-1/2\kappa}.$$ Finally, we can estimate
$c_4$ by
$$c_4(t)\le \|\nabla
{\bf{m}}^\pm_{k,\i,2}\|_{L^\infty(\O)}\|\H_\i\|_{L^\infty([0,T];
L^2(\O))}\|\nabla\H_\i\|_{L^2(\O)}\le C\|\nabla\H_\i\|_{L^2(\O)}.$$
Therefore, using the estimate \eqref{307} repeatedly, we can
conclude that $b^\pm_{k,\i}$ converges strongly to $0$ in
$L^2(0,T)$.

\subsubsection{\bf{The case $k\in J$}}
From \eqref{306} and the fact that $\lambda^\pm_{k,1}=0$, we see
that $\exp^{\pm\textrm{i}\lambda_{k,0}t/\i}b^\pm_{k,\i}$ is bounded
in $L^2(0,T)$ and that its time derivative  is bounded in
$\sqrt{\i}L^1(0,T)+L^p(0,T)$ for some $p>1$. It follows that up to a
subsequence, it converges strongly in $L^2(0,T)$ to some element
$b^\pm_{k,osc}$.

Next, since $\r_\i(x,t)$ converges to $1$ in $C([0,T];
L^\gamma(\O))$ and $b^\pm_{k,\i}(t)$ are uniformly bounded in
$L^2([0,T])$, we deduce that
\begin{equation*}
\r_\i b^\pm_{k,\i}b^\pm_{l,\i}\f{\nabla
\Psi_{k,0}}{\lambda_{k,0}}\otimes\f{\nabla\Psi_{l,0}}{\lambda_{l,0}}-b^\pm_{k,\i}b^\pm_{l,\i}\f{\nabla
\Psi_{k,0}}{\lambda_{k,0}}\otimes\f{\nabla\Psi_{l,0}}{\lambda_{l,0}}\rightarrow
0,
\end{equation*}
in the sense of distributions. Hence, we only need to consider the
terms $$b^\pm_{k,\i}b^\pm_{l,\i}\f{\nabla
\Psi_{k,0}}{\lambda_{k,0}}\otimes\f{\nabla\Psi_{l,0}}{\lambda_{l,0}},\qquad\textrm{for
all } k,l\in J.$$ On the other hand, due to the strong convergence
of $\exp^{\pm\textrm{i}\lambda_{k,0}t/\i}b^\pm_{k,\i}$ in
$L^2([0,T])$ when $k\in J$, we can deduce that
\begin{equation*}
\begin{split}
b^\pm_{k,\i}&b^\pm_{l,\i}\f{\nabla
\Psi_{k,0}}{\lambda_{k,0}}\otimes\f{\Psi_{l,0}}{\lambda_{l,0}}\\&=
\exp^{\textrm{i}(\lambda_{k,0}-\lambda_{l,0})t/\i}\exp^{-\textrm{i}\lambda_{k,0}t/\i}b_{k,0}\exp^{\textrm{i}\lambda_{l,0}t/\i}b_{l,0}
\f{\nabla\Psi_{k,0}}{\lambda_{k,0}}\otimes\f{\nabla\Psi_{l,0}}{\lambda_{l,0}}\\&\rightarrow\exp^{\textrm{i}(\lambda_{k,0}-\lambda_{l,0})t/\i}
b^\pm_{k,osc}(t)b^\pm_{l,osc}(t)\f{\nabla\Psi_{k,0}}{\lambda_{k,0}}\otimes\f{\nabla\Psi_{l,0}}{\lambda_{l,0}},
\end{split}
\end{equation*}
at least in the sense of distributions. Thus, we are only left to
study the interaction of terms
\begin{equation}\label{d1}
\exp^{\textrm{i}(\lambda_{k,0}-\lambda_{l,0})t/\i}
b^\pm_{k,osc}(t)b^\pm_{l,osc}(t)\f{\nabla\Psi_{k,0}}{\lambda_{k,0}}\otimes\f{\nabla\Psi_{l,0}}{\lambda_{l,0}}.
\end{equation}

We will finish the analysis of the interaction by two cases. The
first case is $\lambda_{k,0}=\lambda_{l,0}$. In this case, the term
\eqref{d1} is reduced to
$$b^\pm_{k,osc}(t)b^\pm_{l,osc}(t)\f{\nabla\Psi_{k,0}}{\lambda_{k,0}}\otimes\f{\nabla\Psi_{l,0}}{\lambda_{l,0}},$$
whose divergence is clearly a gradient in the sense of
distributions, due to the fact that as long as
$\lambda_{k,0}=\lambda_{l,0}$, we have
$$\Dv(\nabla\Psi_{k,0}\otimes\nabla\Psi_{l,0}+\nabla\Psi_{l,0}\otimes\nabla\Psi_{k,0})=-\lambda^2_{k,0}\nabla(\Psi_{k,0}\Psi_{l,0})
+\nabla(\nabla\Psi_{k,0}\cdot\nabla\Psi_{l,0}).$$ For the second
case, we have $\lambda_{k,0}\neq\lambda_{l,0}$. Under this
situation, due to the fact that $b^\pm_{k,osc}(t)\in L^2([0,T])$ as
$k\in J$, we know that $b^\pm_{k,osc}(t)b^\pm_{l,osc}(t)\in
L^1([0,T])$ for all $k,l\in J$. Then, by Riemann-Lebesgue Lemma, we
conclude that
$$\int_0^T\exp^{\textrm{i}(\lambda_{k,0}-\lambda_{l,0})t/\i}
b^\pm_{k,osc}(t)b^\pm_{l,osc}(t)dt\rightarrow 0,\quad \textrm{as }
\i\rightarrow 0,$$ which implies that \eqref{d1} converges to $0$ in
the sense of distributions. Hence,  the finite sum, as
$k,l\le M$,
$$\Dv\left(\sum_{k,l\in J}\r_\i b^\pm_{k,\i}b^\pm_{l,\i}\f{\nabla
\Psi_{k,0}}{\lambda_{k,0}}\otimes\f{\nabla\Psi_{l,0}}{\lambda_{l,0}}\right)$$
converges to a gradient in the sense of distributions. And hence,
$$\Dv(\r_\i Q_2\u_\i\otimes Q_2\u_\i)$$ converges to a gradient in
the sense of distributions.


This completes our proof of Theorem \ref{t3}.

\bigskip\bigskip

\section*{Acknowledgments}

Xianpeng Hu's research was supported in part by the National Science
Foundation under grant DMS-0604362. Dehua Wang's research was supported in
part by the National Science Foundation under grant DMS-0604362, and the Office of Naval Research under grant N00014-07-1-0668.

\end{document}